\documentclass[a4paper,12pt]{amsart}
\usepackage[utf8]{inputenc}
\usepackage[english]{babel} 

\usepackage[T1]{fontenc}

\usepackage{amsthm}
\usepackage{amsmath}
\usepackage{amsfonts}
\usepackage{mathtools}
\usepackage[colorlinks]{hyperref}
\usepackage{enumerate}
\usepackage{amssymb}
\usepackage{tikz,pgf}
\usepackage{graphicx}
\usetikzlibrary{matrix,arrows,decorations.pathmorphing}
\usetikzlibrary{calc}
\makeatletter
\newbox\dottedarrow@box
\setbox\dottedarrow@box\hbox
  {%
    \begin{tikzpicture}
      \draw[dotted,->] (0,0) -- (1.5em,0);
    \end{tikzpicture}%
  }
\newcommand*\dottedarrow
  {\relax\ifmmode\expandafter\dottedarrow@m\else\expandafter\dottedarrow@t\fi}
\newcommand*\dottedarrow@t[1][1.5em]
  {\resizebox{#1}{!}{\raisebox{.5ex}{\usebox\dottedarrow@box}}}
\newcommand*\dottedarrow@m[1][]
  {%
    \if\relax\detokenize{#1}\relax
      \mathchoice
        {\dottedarrow@t}
        {\dottedarrow@t}
        {\dottedarrow@t[1.1em]}
        {\dottedarrow@t[0.9em]}%
    \else
      \dottedarrow@t[#1]%
    \fi
  }
\makeatother

\def\bdi{\begin{diagram}}
\def\edi{\end{diagram}}

\newtheorem{thm}{Theorem}[section]
\newtheorem{cor}[thm]{Corollary}
\newtheorem{lem}[thm]{Lemma}
\newtheorem{prop}[thm]{Proposition}
\theoremstyle{definition}
\newtheorem{defi}[thm]{Definition}
\newtheorem{defis}[thm]{Definitions}
\newtheorem{conj}[thm]{Conjecture}
\newtheorem{conv}[thm]{Convention}
\newtheorem{nota}[thm]{Notation}
\newtheorem{rem}[thm]{Remark}
\newtheorem{rems}[thm]{Remarks}
\newtheorem{exa}[thm]{Example}
\newtheorem{exas}[thm]{Examples}
\newtheorem{prob}[thm]{Problem}
\newtheorem{probs}[thm]{Problems}
\newtheorem{ques}[thm]{Question}
\newtheorem{sett}[thm]{Setting}
\newtheorem{sit}[thm]{}

\newcommand{\Hol}{ \operatorname{{\rm Hol}}}

\newcommand{\Aut}{ \operatorname{{\rm Aut}}}

\def\deg{\mathop{\rm deg}}

\renewcommand{\epsilon}{\varepsilon}

\def\and{\quad\mbox{and}\quad}

\newcommand{\C}{\ensuremath{\mathbb{C}}}

\newcommand{\R}{\ensuremath{\mathbb{R}}}

\newcommand{\G}{\ensuremath{\mathbb{G}}}

\newcommand{\bX}{{\bar X}}

\newcommand{\bC}{{\bar C}}

\newcommand{\bY}{{\bar Y}}

\newcommand{\tL}{{\tilde L}}

\newcommand{\tE}{{\tilde E}}

\newcommand{\bL}{{\bar L}}

\newcommand{\cL}{{\ensuremath{\mathcal{L}}}}

\newcommand{\cE}{{\ensuremath{\mathcal{E}}}}

\renewcommand{\rho}{\varrho}

\def\bals#1\eals{\begin{align*}#1\end{align*}}
\def\bal#1\eal{\begin{align}#1\end{align}}

\def\SAut{\mathop{\rm SAut}}

\def\PP{{\mathbb P}}

\def\dd{ {\rm d }}

\renewcommand{\phi}{\varphi}

\newcommand{\bnum}{\begin{enumerate}}
\newcommand{\enum}{\end{enumerate}}

\newcommand{\brem}{\begin{rem}}
\newcommand{\brems}{\begin{rems}}
\newcommand{\erem}{\end{rem}}
\newcommand{\erems}{\end{rems}}
\newcommand{\bprob}{\begin{prob}}
\newcommand{\eprob}{\end{prob}}
\newcommand{\bprobs}{\begin{probs}}
\newcommand{\eprobs}{\end{probs}}
\newcommand{\bques}{\begin{ques}}
\newcommand{\eques}{\end{ques}}
\newcommand{\bexa}{\begin{exa}}
\newcommand{\bexas}{\begin{exas}}
\newcommand{\eexa}{\end{exa}}
\newcommand{\eexas}{\end{exas}}
\newcommand{\bdefi}{\begin{defi}}
\newcommand{\edefi}{\end{defi}}
\newcommand{\bdefis}{\begin{defis}}
\newcommand{\edefis}{\end{defis}}
\newcommand{\bcor}{\begin{cor}}
\newcommand{\ecor}{\end{cor}}
\newcommand{\blem}{\begin{lem}}
\newcommand{\elem}{\end{lem}}
\newcommand{\bconv}{\begin{conv}}
\newcommand{\econv}{\end{conv}}
\newcommand{\bconj}{\begin{conj}}
\newcommand{\econj}{\end{conj}}
\newcommand{\bprop}{\begin{prop}}
\newcommand{\eprop}{\end{prop}}
\newcommand{\bthm}{\begin{thm}}
\newcommand{\ethm}{\end{thm}}
\newcommand{\bnota}{\begin{nota}}
\newcommand{\enota}{\end{nota}}
\newcommand{\bsit}{\begin{sit}}
\newcommand{\esit}{\end{sit}}
\newcommand{\be}{\begin{equation}}
\newcommand{\ee}{\end{equation}}
\newcommand{\bproof}{\begin{proof}}
\newcommand{\eproof}{\end{proof}}
\newcommand{\bsett}{\begin{sett}}
\newcommand{\esett}{\end{sett}}
\def\ba{\begin{array}}
\def\ea{\end{array}}

\begin{document}

\title[Kobayashi-Royden pseudometric]{Kobayashi-Royden pseudometric of contractible surfaces}

\author{Shulim Kaliman}

\address{Department of Mathematics,
University of Miami, Coral Gables, FL 33124, USA}
\email{kaliman@math.miami.edu}

\dedicatory{To Leonid Makar-Limanov on the occasion of his 80-th birthday}

\thanks{2020 \emph{Mathematics Subject Classification} 
Primary  14} 

\date{}
\maketitle
\begin{abstract}  We describe a family of smooth contractible algebraic surfaces $X$ different from $\C^2$
such that $X$ admits  dominant holomorphic maps from $\C^2$ and there is a unique line $E$ in $X$
for which the Kobayashi-Royden pseudometric vanishes on the tangent bundle over $X\setminus E$.
\end{abstract}

\section{Introduction}

All varieties in this paper are considered over the field of complex numbers.
In  1971 Ramanujam \cite{Ra71}
 constructed the first smooth contractible affine algebraic surface different from $\C^2$.
Such surfaces $X$ are often called Ramanujam surfaces \cite{Za90}. 
 Fujita \cite{Fu79} and Miyanishi and Sugie \cite{MiSu80}
 showed that $X\times \C^n$ is not isomorphic
to $\C^{n+2}$ for every $n>0$, thus, proving the Zariski cancellation conjecture in
dimension 2. Later Fujita \cite{Fu82} developed an effective  method of constructing Ramanujam surfaces.
But his construction (as well as the Ramanujam original paper) produces only
 Ramanujam surfaces $X$ with logarithmic Kodaira 
dimension $\bar \kappa (X)= 2$.  Gurjar and Miyanishi \cite{GuMi88} discovered 
Ramanujam surfaces $X$ with $\bar \kappa (X)= 1$
 and found their classification. They also showed that there is no Ramanujam surface $X$
with $\bar \kappa (X)= -\infty$, whereas the result of Pe{\l}ka and Ra\'zny \cite{PeRa21} implies
that $X$ cannot have $\bar \kappa (X)= 0$ (see Proposition \ref{pre.p1}).
An elegant description of Ramanujam surfaces $X$ with $\bar \kappa (X)= 1$ was
found later by Petrie and tom Dieck \cite{tDiPe90} (actually, all of them can be presented as explicit hypersurfaces 
in $\C^3$ \cite{KaML96}).

On the other hand, the result of Choudary and Dimca \cite{ChDi94} implies that 
$X\times \C^n$ is diffeomorphic to $\R^{2n+4}$ for $n\geq 1$.
So, it is natural to ask whether $X\times \C^n$ is biholomorphic to $\C^{n+2}$.
A partial answer to this question is based on the paper of
Sakai \cite{Sa77} who proved that every complex 
algebraic variety of general type is measure hyperbolic (that is, its top Eisenman measure is almost
nowhere zero).  Using this fact Zaidenberg \cite{Za90}, \cite{Za93} showed that 
 $X\times \C^n$ is not biholomorphic to $\C^{n+2}$ when $\bar \kappa (X)= 2$.
 However, Ramanujam surfaces $X$ with $\bar \kappa (X)= 1$ are not measure hyperbolic
 and the  question whether  $X\times \C^n$ is not biholomorphic to $\C^{n+2}$ remains open. 
 The answer to this question would be positive is $X$ possesses a nontrivial Kobayashi-Royden pseudometric which is also the first Eisenman measure \cite{Ka94}
 (for the definition of Eisenman measures we refer to \cite{Ko76} and \cite{Ei70}).
 In particular, one may hope that as in the case of the top Eisenman measures for varieties of general type
 the Kobayashi-Royden pseudometric is  almost nowhere zero on tangent bundle $TX$.
 An unsuccessful attempt to prove this fact was the original motivation for our paper.
 Actually, we find a family of Ramanujam surfaces $X$
with $\bar \kappa (X)=1$ and Kobayashi-Royden pseudometric vanishing on a dense
open subset of $TX$ (see Section \ref{s.kob}).  Furthermore, we show  that 
there a holomorphic map $\C^2 \to X$ such that its Jacobian determinant
 is not identically zero (a surface admitting such map is called
 holomorphically dominable by $\C^2$),  
 whereas there is no dominant morphism
from $\C^2$ to $X$.
The surfaces  holomorphically dominable
 by $\C^2$ were extensively studied by Buzzard and Lu \cite{BuLu00}. 
However, the families we discovered are not in their list.

\section{Preliminaries} 

Throughout this paper a notation of the
form $\C_{z_1, \ldots, z_n}^n$ means $\C^n$ equipped with a coordinate system $(z_1, \ldots, z_n)$ 
(in particular, $\C_z$ is a line equipped with a coordinate $z$). 
As usual, for a complex manifold $X$ its tangent bundle is denoted by $TX$ and $T_xX$ is the tangent space 
at $x \in X$.

\bdefi\label{}  For $r>0$ let $\Delta_r \subset \C_z$ be the disc of radius $r$ with the center at the origin $0$
and $\Hol (\Delta_r, X)$ be the set of holomorphic maps from $\Delta_r$ to $X$.
Given a vector $\nu \in T_xX$ its Kobayashi-Royden pseudometric  is
$$K_X (\nu )= \inf \{ 1/r \, |\, \phi \in \Hol (\Delta_r, X) \text{ with } \phi (0)=x \text{ and } \dd \phi \big(\frac{\dd}{\dd z}\big|_0\big)=\nu \}.$$
\edefi

In particular, if one has a holomorphic map $\phi : \C \to X$ with $\phi (0)=x$ and  $\phi \big(\frac{\dd}{\dd z}\big|_0\big)=\nu$,
then $K_X (\nu )=0$.

Further in this paper every complex manifold $X$ will be the complement to
an effective divisor $D=D_1+ \ldots + D_s$ of
simple normal crossing (SNC) type in a projective 
algebraic manifold $\bX$. Let $n=\dim X$ and $\Omega_n(X)$ be the sheaf of germs
of holomorphic $n$-forms on $X$. In particular, the space $\Gamma (X,\Omega_n(X))$ of holomorphic sections of $\Omega_n(X)$ 
is the space of top holomorphic differential forms.
 For every point
of $\bX$ there exists a local coordinate system
$(z,w)=(z_1, \ldots , z_l, w_1, \ldots , w_{n-l})$
such that 
$$z_1 \cdots z_l=0$$
defines the germ of $D$
around this point ($0\leq l \leq n$ and
when $l=0$ this point does not belong
to $D$).
In this case $\Omega_n (X)$
contains the subsheaf $\cL=\Omega_n (\bar {X}, D)$ of logarithmic $n$-forms
along $D$.
The germs of these forms can be written as
$$ \sum_{r+q=k} c_{I,J} (z,w) {\frac {dz_{i(1)}}{z_{i(1)}}} \wedge
\ldots {\frac {dz_{i(r)}}{z_{i(r)}}} \wedge
{\frac {dw_{j(1)}}{w_{j(1)}}} \wedge
\ldots {\frac {dw_{j(q)}}{w_{j(q)}}},$$
where $c_{I,J}(z,w)$ is the germ of a holomorphic function, 
and the indices of summation are
$I=(i(1), \ldots , i(r))$ and
$J=(j(1), \ldots , j(q))$. 

\bdefi
If for every $m>0$ the $m$th power $\cL^{m}$ of $\cL$ does not have a non-trivial global
section, then the logarithmic Kodaira dimension of $X$ is $\bar \kappa (X) = - \infty$. Otherwise,
$$\bar \kappa (X)= \lim_{m \to + \infty } \sup {\frac {\log \dim \Gamma (X, \cL^{m})}
{\log m}} .$$
 \edefi

The logarithmic Kodaira dimension does not depend of the choice of SNC completion $\bX$ of $X$
and one also has $\bar \kappa (X) \leq n$  (\cite{Ii77}).

Every smooth contractible affine algebraic surface $X$ is factorial \cite{Fu82} and
the ring $\C [X]^*$ of invertible functions on it consists of constants only.
Gurjar and Miyanishi \cite{GuMi88} classified all smooth factorial surfaces $X$ with $\C [X]^*=\C^*$
and logarithmic Kodaira dimension equal to $\bar \kappa (X)=1$ 
(this yields, in particular a classification of Ramanujam surfaces with logarithmic Kodaira dimension 1,
whereas for negative dimension such Ramanujam surfaces do not exist). However, their classification for the case 
$\bar \kappa (X)=0$ contains a flaw. Freudenburg, Kojima, and Nagamine \cite{FrKoNa19} found an
 infinite collection of smooth affine factorial surfaces $X$ with  $\C [X]^*=\C^*$ and $\bar \kappa (X)=0$
that are not in the list presented in \cite{GuMi88}. Later Pe{\l}ka and Ra\'zny \cite{PeRa21} obtained a complete classification
of such surfaces. The next proposition is an  immediate consequence of  the latter result.

\bprop\label{pre.p1} There are no Ramanujam surfaces $X$ with $\bar \kappa (X)=0$.
\eprop

\bproof Let $X$ be a smooth affine factorial surface with  $\C [X]^*=\C^*$ and $\bar \kappa (X)=0$.
By \cite[Theorem 1.1]{PeRa21} there exist  monic polynomials $p_1(t)$ and $p_2(t)$ such that $X$ is isomorphic
to the spectrum of the ring 
$$\C [x_1,x_2][ (x_2x_1^{-\deg p_1}-p_1(x_1^{-1}))x_1^{-1},  (x_1x_2^{-\deg p_2}-p_2(x_2^{-1}))x_2^{-1}  ].$$
This implies that the restriction of the natural morphism $\phi : X\to \C^2_{x_1,x_2}$ is isomorphism over
the torus $T=\{ (x_1,x_2) \in \C^2| \, x_1, x_2 \ne 0\}$, the image of $\phi$ is $T\cup \{ (1,0), (0,1)\}$, 
and $\phi^{-1} ((1,0))\simeq \phi^{-1} ((0,1))\simeq \C$. Thus, $X$ can be stratified into a disjoint union
of $T$ and two lines. Since the Euler characteristic of $T$ is 0 and of $\C$ is 1, the Euler characteristic of $X$ is 2 by \cite{Du87}
which yields the desired conclusion.
\eproof

\section{Ramification}\label{s.ram}

Let $L$ be the curve in $\C^2_{u,v}$ given by $u^k-v^l=1$ where $k$ and $l$ are fixed relatively prime natural 
numbers with $k>l \geq 2$ 
and $L_0\subset \C^2_{u,v}$ be given by $u^k-v^l=0$.
The following fact can be extracted from \cite[p. 150-151]{tDiPe90}.

\bprop\label{ram.p1} {\rm (tom Dieck-Petrie)}.
Let $\rho_1 : V_{1}' \to \C_{u,v}^2$ be the monoidal transformation at  the point $(1,1)$, $E_1'$ be its exceptional divisor,
$\tL_0 \subset V_{1}' $ be the proper transform of $L_0$, $V_{1}=V_{1}' \setminus \tL_0$, and
$E_1=E_1'\setminus \tL_0$. Suppose that $V_{i}' \to V_{i-1}', \, i=2, \ldots, m$ is a monoidal transformation at 
a point of $E_{i-1}$, $E_i'$ is its exceptional divisor, $\tE_{i-1}\subset V_{i}' $ is the proper transform of 
$E_{i-1}'$, $V_{} =V_{i}' \setminus \tE_{i-1}$, and $E_i=E_i'\setminus \tE_{i-1}$.
Then every $V_{i}, \, i=1, \ldots, m$ is a Ramanujam surface with $\bar \kappa (V_{i})=1$.
Furthermore, up to an isomorphism every Ramanujam surface with logarithmic Kodaira dimension 1
can be obtained via such construction.
\eprop

\blem\label{ram.l1}  Let $\rho_m : V_{m} \to \C_{u,v}^2$ be the natural morphism, that is,
$V_{m}$ is disjoint union of  $\C^2\setminus L_0$ and $E_m\simeq \C$.
Then there are $c_1, \ldots , c_{m-1} \in \C$ such that the  function 
$\rho_m^* (\frac{u-1-\sum_{i=1}^{m-1}c_i(u^k-v^l)^i}{(u^k-v^l)^{m}})$ 
is regular on $V_{m} $ and its restriction yields a coordinate on $E_m$.
\elem

\bproof 
Let $I_0$ be the defining ideal of the point $p_0=(1,1) \in \C^2$.
Note that the localization of $I_0$ at $p_0$ is generated by $u-1$ and $u^k-v^l$. 
Since $L_0=(u^k-v^l)^*(0)$, the function $\frac{u-1}{u^k-v^l}$ is regular on $V_{1}$ and
 its restriction yields a coordinate on $E_1$. Note that  the lift of $u^k-v^l$ vanishes on $E_1$ with multiplicity 1.
  To get $V_{2}$ one needs to blow up $V_{1}$  at a point $p_1$ in $E_1$ with a coordinate $c_1$. 
  Hence, the localization of the ideal  of the blowing-up  at $p_1$
  is generated by $\frac{u-1}{u^k-v^l}-c_1$ and $u^k-v^l$. This implies that
  the restriction of $\frac{u-1-c_1(u^k-v^l)}{(u^k-v^l)^2}$ yields a coordinate on $E_2$.
  Proceeding in this manner by induction one gets the general case. 
\eproof

The following fact is straightforward.

\blem\label{ram.l2}
Let $Y= L  \times \C_z$.
Consider the morphism $\chi : Y \to \C^2_{x,y}, \, (u,v,z) \mapsto (uz^l, vz^k)=:(x,y)$ and $Y^*=\C^*\times L$.
Then $\chi(Y^*)$ is the complement to $L_0=(x^k-y^l)^*(0)\subset \C_{x,y}^2$ and $\chi(Y^*)$ is isomorphic to the quotient $Y^*/\mu_{kl}$ where 
$\mu_{kl}$ is the multiplicative group generated by a primitive $kl$-root of unity
and the action of $\mu_{kl}$ on $Y$ is given by $\lambda . (u,v,z)=(\lambda^lu, \lambda^kv, \lambda^{-1} z)$. 
\elem
 
 \blem\label{ram.l3}  Let $\bY_0$ be a smooth surface equipped with a $\mu_{kl}$-action,
 $C_0\simeq \C$ be a curve in $\bY_0$, $z$ and $\zeta$ be functions regular in a neighborhood
 of $C_0$ in $\bY_0$ such that $z$ vanishes on $C_0$ with multiplicity 1 and the restriction
 of $\zeta$ is a coordinate on $C_0$. Suppose that for every $\lambda \in \mu_{kl}$
 one has $\lambda.\zeta=\zeta$ and $\lambda.z =\lambda^{-1}z$.
  Let $\delta: \bY_1 \to \bY_0$ be
a sequence of $kl$ monoidal transformations  of $\bY_0$ at the point in $C_0$ 
with coordinate $\zeta=1$  and infinitely near points, $\bC_i$ be the curve generated by  the $i$th monoidal transformation,
and the center of  the  $(i+1)$st  monoidal transformation occurs at a point of
$C_i=\bC_i\setminus \bC_{i-1}$ (where we identify $\bC_i$ with its proper transform). Suppose that
the $\mu_{kl}$-action induces a regular action on $\bY_1$. Then
the lift of $z$ vanishes on $\bC_{kl}$ with multiplicity 1,
the function  $\frac{\zeta-1}{z^{kl}}$ is regular in a
neighborhood of $C_{kl}$, and its restriction yields a coordinate on $C_{kl}\simeq \C$.
\elem

\bproof After the first monoidal transformation $\frac{\zeta-1}{z}$ is the coordinate function on $C_1\simeq \C$
and the lift  of $z$ is a function that vanishes on $C_1$ with multiplicity 1. 
Note that $\frac{\zeta-1}{z}$ is not $\mu_{kl}$-invariant. 
Hence, the  second monoidal must occur at the point of $C_1$ with $\frac{\zeta-1}{z}=0$ since any other point 
of $C_1$ is not preserved by the induced $\mu_{kl}$-action on $C_1$.
Then $\frac{\zeta-1}{z^2}$ is a coordinate function on $C_2$ and
the lift  of $z$ is a function that vanishes on $C_2$ with multiplicity 1. 
Note that $\frac{\zeta-1}{z^2}$ is not $\mu_{kl}$-invariant. Thus, by induction we see that after
$j$ monoidal transformations $\frac{\zeta-1}{z^j}$ is a coordinate function on $C_j$
and the lift  of $z$ is a function that vanishes on $\bC_j$ with multiplicity 1. 
Furthermore,  $\frac{\zeta-1}{z^j}$ is not $\mu_{kl}$-invariant unless $j=kl$ which concludes the proof.
\eproof

\bnota\label{ram.n1}
Let $Y$ be as in Lemma \ref{ram.l2} and $\bL$ be the completion of $L$ smooth at $o =\bL \setminus L$.
Then one can suppose that a local coordinate $w$ on $\bL$ at  $o$ is such that $u=w^{-l}$.
Note that the $\mu_{kl}$-action sends $w$ to $\lambda^{-1} w$. 
Furthermore, $Y$ admits a natural embedding into  $\bY=\bL \times  \PP$ with the boundary divisor $D = \bY \setminus Y$
containing the irreducible component $o \times \PP^1$ given locally by $w=0$. 
Note that the $\mu_{kl}$-action extends to an action on $\bY$.
Consider the monoidal transformation $\bY_0\to \bY$ at the point with coordinates $z=w=0$. 
Then the exceptional curve $\bC_0$ contains an open subset $C_0\simeq \C$ such that  $\zeta = \frac{z}{w}$
is regular in a neighborhood of $C_0$ in $\bY_0$ and its restriction yields a coordinate on $C_0$.
Observe that  $\zeta$ is invariant under the $\mu_{kl}$-action.
Hence, the $\mu_{kl}$-action on $\bY$ induces an action on $\bY_0$ such that its restriction to $C_0$
is the trivial action.
\enota

\blem\label{ram.l4}  
 Let $\delta_m: \bY_m \to \bY_0$ be
a sequence of $mkl$ monoidal transformations  of $\bY_0$ at the point in $C_0$ 
with coordinate $\zeta=1$  and infinitely near points,
 $\bC_i$ be the curve generated by  the $i$th monoidal transformation,
and the center of  the  $(i+1)$st  monoidal transformation occurs at a point of
$C_i=\bC_i\setminus \bC_{i-1}$ (where we identify $\bC_i$ with its proper transform). 
 Suppose that  the $\mu_{kl}$-action induces an action on $\bY_m$. Then for some $c_1, \ldots , c_{m-1} \in \C$
the function  $\frac{\zeta-1-\sum_{i=1}^{m-1}c_iz^{ikl}}{z^{mkl}}$ is regular in a
neighborhood of $C_{mkl}$, and its restriction yields a coordinate on $C_{mkl}\simeq \C$.
\elem

\bproof  Lemma \ref{ram.l3} yields the case $m=1$.  In particular,  $\frac{\zeta-1}{z^{kl}}$ is regular in 
a neighborhood of $C_{kl}$, its restriction is a coordinate on $C_{kl}$,
and the lift of $z$ vanishes on $C_{kl}$ with multiplicity 1.
For $m=2$ we observe that  the $(kl+1)$st monoidal transformation
 cannot occur at a point of $C_i\setminus \bC_{i+1}$ for $1\leq i\leq kl-1$ since 
such point is not preserved by the $\mu_{kl}$-action.
Therefore,  this monoidal transformation  must
occur at some point of $C_{kl}$ with a coordinate $\frac{\zeta-1}{z^{kl}}=c_1$. Hence, by Lemma \ref{ram.l3}
 the restriction of $\frac{\zeta-1-c_1z^{kl}}{z^{2kl}}$ yields a coordinate on $C_2$.
 Proceeding in this manner by induction one gets the general case. 
\eproof

\bprop\label{ram.p2} Let $Y_m$
be obtained from $\bY_m$ by removing the proper transform  $D_m$ of $D$ and the curves $C_0, \ldots , C_{mkl-1}$,
so, $Y_m$ is the disjoint union of $\delta_m^{-1}(Y^*)\simeq Y^*$ and $C_{mkl}$.     
Let $V_{m}=(\C^2\setminus L_0)\cup E_m$ be as in Lemma \ref{ram.l1}.
Then the rational map $\chi_m: Y_m\dashrightarrow V_{m}$ induced by $\chi: Y \to \C^2$ is regular and it maps
$C_{mkl}$ isomorphically onto $E_m$.
\eprop

\bproof  By Lemma \ref{ram.l1}  the lift of $\frac{x-1-\sum_{i=1}^{m-1}c_i(x^k-y^l)^i}{(x^k-y^l)^{m}}$ is regular
on $V_{m}$ and its restriction yields a coordinate on $E_m$.
The pullback of this function to $Y_m$ is $\frac{\zeta-1-\sum_{i=1}^{m-1}c_iz^{ikl}}{z^{mkl}}$.
By Lemma \ref{ram.l4} the latter function is regular in a
neighborhood of $C_{mkl}$, and its restriction yields a coordinate on $C_{mkl}$.
Since by construction $\chi_m$ is regular over $\C^2\setminus L_0$ this yields the desired conclusion.
\eproof

\section{Kobayashi-Royden pseudometric}\label{s.kob}

In this section use the notation as in Section \ref{s.ram} with  $k=3$ and $l=2$. In particular, $\bL$ is a smooth elliptic curve.

\bthm\label{kob.t1}  Let $V_m$ be as in Proposition \ref{ram.p1} with $k=3$ and $l=2$.
Let $V_m^*=V_{m} \setminus E_m$ and $TV_m^*$ be the tangent bundle on $V_m^*$.
Then  for every $\nu \in TV_m^*$ there is a holomorphic
map $h: \C \to V_{m}$ such that $h_*(\frac{\dd}{\dd \xi}\big|_0)=\nu$ where $\xi$ is a coordinate on $\C$.
\ethm

We start the proof with the following lemma.

\blem\label{kob.l1} Let $\Gamma$ be the lattice of Gaussian integers in $\C_\xi$
and $P(\xi)$ be a polynomial of degree $n$ such that $P(0)=0$ and $P'(0)=1$. 
Let $p\in \C \setminus \Gamma$.
Then there is an entire function
$h(\xi)$ such that 

{\rm  (i)} $h$ vanishes on $\Gamma$ only;

{\rm  (ii)} for every $q \in \Gamma$ the first $n+1$ terms of the Taylor expansion of $h$ at $q$
coincide with $P(\xi-q)$;

{\rm  (iii)} the Taylor expansion of $h$ at $p$ coincides with a prescribed polynomial.
\elem

\bproof The Weierstrass factorization theorem implies that 
$h(\xi)=e^{g(\xi)}\xi\prod_{q \in \Gamma\setminus \{ 0\}} \cE_q	 (\xi)$
where each factor $\cE_q$ vanishes at $	$ only and its derivative at $q$ is nonzero.
Hence, choosing the entire function $g(\xi)$ with appropriate Taylor expansions at the points of $\Gamma \cup \{ p \}$ 
we get the desired conclusion.
\eproof

\blem\label{kob.l2} Let $\pi : \C \to \bL$ be a universal covering and $\Gamma=\pi^{-1}(o)$ (where $o=\bL\setminus L$).
Let $\phi : \C_{\xi}\to \bY=\bL\times \C, \, \xi \mapsto (\pi(\xi), h(\xi))$ and $\sigma_m: Y_m \to \bY$ be the natural projection
where $Y_m$ is as in Proposition \ref{ram.p2}.
Then for an appropriate choice of $P$ from Lemma \ref{kob.l1} there exists $\psi : \C \to Y_m$ such that 
$\phi =\sigma_m \circ \psi$.
\elem

\bproof  By Lemma \ref{kob.l1} $\phi^{-1} (\bL \times 0) =\Gamma$ and $\phi (p)=(o, 0)$ for every $p \in \Gamma$.
Furthermore, for $\bY_0$, $C_0$, and $\zeta$ as in Notation \ref{ram.n1}
there is a lift $\chi : \C \to \bY_0$  of $\phi$ such that $\chi (p)$
is a point in $C_0$ and the coordinate  of this point is $\zeta =1$ since $P'(p)=1$. Hence, choosing appropriate $P$
one can lift $\chi$ to a holomorphic map $\psi : \C \to Y_m$ which concludes the proof.
\eproof

\proof[Proof of Theorem \ref{kob.t1}]  Let $Y_m^*=Y_m\setminus C_{6m}$. 
Since the morphism $\chi_m: Y_m\to V_{m}$ is smooth finite over $V_m^*$ it suffices to show that for
a given $\nu \in TY_m^*$ there is a holomorphic map $\psi : \C_\xi \to Y_m$ such that 
$\psi_*(\frac{\dd}{\dd \xi}\big|_0)=\nu$. Let $\psi$ and $\phi=\sigma_m\circ \psi$ be as in Lemma \ref{kob.l2}
and, in particular, $\phi=(\pi, h)$. 
We can suppose that $\psi (0)$ coincides with the image $q\in Y_m$ of $\nu$ under the natural projection $TY_m\to Y_m$.
We use the natural identification of  $T_qY_m$ with $\C^2$ induced by the isomorphism $Y_m^*\simeq L\times \C^*$. 
Since  the restriction of $\sigma_m$ over $Y_m^*$ can be viewed as the identity map and $\pi$ is smooth
we have $\psi_*(\frac{\dd}{\dd \xi}\big|_0)=(\nu_1, \nu_2) \in T_qY_m$
where $\nu_1 \ne 0$.  By Lemma \ref{kob.l1}(iii) $\nu_2$ can be an arbitrary number. This yields the claim for
vectors $\nu=(\nu_1, \nu_2)$ with $\nu_1 \ne 0$. For $\nu_1=0$ the claim follows from the fact that there is
a $\C^*$-curve tangent to $\nu$ since
$Y_m^*\simeq \bL\times \C^*$. Hence, we are done. \hfill $\square$

\bcor\label{c1} The Kobayashi-Royden pseudometric of $V_m$ vanishes on  $TV_{m}^*$.

\ecor

\brem\label{kob.r1} Let $g_0$ be a nonzero function vanishing at every point of $\Gamma \cup \{ p \}$  with sufficiently high degree.
Note that every function of the form $g_c=g +c g_0, \, c \in \C$ can be used instead of $g$ in the proof of Lemma \ref{kob.l1}.
In particular, we can make the function $h$ in Lemma \ref{kob.l1} depending on $c$ by letting
$h_c(\xi)=e^{g_c(\xi)}\xi\prod_{p \in \Gamma\setminus \{ 0\}} \cE_p (\xi)$. This in turn yields in  Lemma \ref{kob.l2} the morphisms
$\phi_c : \C_{\xi}\to \bY=\bL\times \C, \, \xi \mapsto (\pi(\xi), h_c(\xi))$ and
$\psi_c : \C \to Y_m$ such that  $\phi_c =\sigma_m \circ \psi_c$ depending holomorphically on $c$.

\erem


\blem\label{kob.l3} The image of the holomorphic map $\Psi :\C_{\xi,c}^2\to Y_m,  (\xi, c) \mapsto \psi_c(\xi)$
is dense in the standard topology.
\elem

\bproof Since $\phi_c =\sigma_m \circ \psi_c$  it suffices to show that the map
$\Phi :\C_{\xi,c}^2\to  \bY,  (\xi, c) \mapsto \phi_c(\xi)= (\pi(\xi), h_c(\xi))$ has the image $I$
dense in the standard topology. Note that $h_c(\xi) =e^{cg_0(\xi)}h_0(\xi)$.
Hence, for every $\xi_0 \in \C_\xi \setminus \Gamma$ such that $g_0(\xi_0)\ne 0$ the image $I$
contains $h_0(\xi_0) \times \C^*\subset \bL\times \C^*=Y_m^*$. This yields the desired conclusion.
\eproof

\bprop\label{kob.p2} Let $V_m$ be as in Theorem \ref{kob.t1} and, so, 
$V_m\setminus E_m \simeq \C^2_{u,v}\setminus \{ (u,v)| \, u^3-v^2=0\}$.
Then $V_m$ is  holomorphically dominable by $\C^2$.
\eprop

\bproof  By Proposition \ref{ram.p2} we have the surjective morphism $\chi_m : Y_m \to V_m$.
Hence,  Lemma \ref{kob.l3} implies the desired conclusion.
\eproof

\brem\label{kob.r2} 
(1) It is interesting to compare the abundance of holomorphic maps $\C \to V_m$ with 
nonconstant morphisms from $\C$ to $V_m$. Gurjar and Miyanishi \cite{GuMi88} 
proved that $E_m$ is the only line contained in $V_m$. Furthermore,
the image of every nonconstant morphism
$\C\to V_{m}$ is contained in $E_m$ \cite[Corollary 3.3]{KaML96}.

(2) In particular, $V_m$ does not admit a dominant morphism from $\C^2$.
It is easy to find rich families of affine surfaces holomorphically dominable by $\C^2$
that do not admit a dominant morphism from $\C^2$.
Say, generalized Gizatullin  surfaces \cite{KaKuLe20} belong to this class.
The most extensive study on this subject is due to Buzzard and Lu \cite{BuLu00}.
In particular, they developed several sufficient criteria for an affine surface to be holomorphically domianble by $\C^2$.
For affine surfaces with logarithmic Kodaira dimension 1 the corresponding criterion is described in \cite[Theorem 5.10]{BuLu00}.
However, the family in Proposition \ref{kob.p2} does not satisfy the assumptions of that theorem.

(3) Let $e: \C \to \C$ be an entire nonconstant function, $\tilde \pi=\pi\circ e$, and $\tilde \Gamma =e^{-1} (\Gamma)$.
Then one can replace $\Gamma$ in the formulation of Lemma \ref{kob.l1} by $\tilde \Gamma$.
Consequently $\phi$ in Lemma \ref{kob.l2}  can by replaced
by $\tilde \phi =(\tilde \pi, h)$ and we get $\tilde \psi : \C \to Y_m$ such that 
$\tilde \phi =\sigma_m \circ \tilde \psi$.  Actually, every holomorphic map $\C\to Y_m$ whose image
is different from $q\times \C^*, \, q \in \bL$ is of this form.
\erem

\bprop\label{kob.p3} Let $\chi_m: Y_m \to V_m$ be as in Proposition \ref{ram.p2} and $\tilde \psi : \C \to Y_m$ be
as in Remark \ref{kob.r2} (3). Let $C=\chi_m\circ \tilde \psi (\C)$. The the order of tangency between $C$ and $E_m$
at all points of $C\cap E_m$ is at least $6$.
\eprop

\bproof This follows from the fact that $\chi_m|_{C_{6m}} : C_{6m} \to E_m$ is an isomorphism and 
$\chi_m : Y_m \to V_m$ is a finite morphism ramified along $E_m$ with ramification index $kl=6$.
\eproof

This leads to the question whether the Kobayashi-Royden pseudometric of any vector from $TV_m|_{E_m} \setminus TE_m$ is nontrivial. \\

{\em Acknowledgement.} It is a pleasure to thank M. Zaidenberg who informed me about the flaw in the paper \cite{GuMi88}
and drew my attention to the papers \cite{FrKoNa19} and \cite{PeRa21}.

\end{document}